\documentclass[a4paper,11pt,twoside]{article} 
\usepackage{amsmath}
 \usepackage{amsthm} 
\usepackage{amssymb}
\usepackage{amsfonts}
\usepackage{latexsym}
\usepackage{makeidx}
\usepackage{mathrsfs}
\usepackage{vruler}
\usepackage{amsxtra}
\usepackage[bf,small]{caption}
\usepackage{sidecap}

\usepackage{graphics}
\usepackage{graphicx}
\usepackage{color,rotating}
\usepackage{epsfig}

\let\oldcaption\caption  
\renewcommand{\caption}[1]{\oldcaption{\bf #1}}

\newtheorem{definition}{Definition}[section]

\newtheorem*{lemma*}{Lemma}

\newtheorem{theorem}[definition]{Theorem}

\newtheorem*{problem*}{Problem (P)}
\newtheorem*{forbiddenballs*}{Forbidden Geodesic Balls (FGB)}
\newtheorem*{touching*}{Touching Principle (TP)}

\newtheorem*{touchingsituations*}{Possible local touching situations}

\newtheorem*{c1*}{(C1)}

\newtheorem*{c2*}{(C2)}

\newcounter{detail}[definition]

\newcommand{\dist}{{\rm dist}}
\newcommand{\area}{{\rm area}}
\newcommand{\distS}{{\dist_{\mS^2}}}

\newfont{\thickmath}{msbm10 scaled \magstephalf}%
\newfont{\smallthickmath}{msbm7 scaled \magstephalf}%
\newfont{\footnotethickmath}{msbm8}%
\newfont{\footnotesmallthickmath}{msbm6}%

\def\bbbs{\mbox{\thickmath S}}

\def\bbbr{\mbox{\thickmath R}}
\def\R{\mbox{\thickmath R}}
\def\bbbz{\mbox{\thickmath Z}}
\def\bbbn{\mbox{\thickmath N}}

\def\mbbbs{\mbox{\smallthickmath S}}

\renewcommand{\S}{\bbbs}
\newcommand{\Z}{\bbbz}
\newcommand{\mS}{\mbbbs}

\renewcommand{\R}{\bbbr}

\newcommand{\N}{\bbbn}

\def\cB{\mathscr{B}}
\def\cT{\mathscr{T}}

\def\mvint_#1{\mathchoice
          {\mathop{\vrule width 6pt height 3 pt depth -2.5pt
               \kern -8pt \intop}\nolimits_{\kern -3pt #1}}%
    {\mathop{\vrule width 5pt height 3 pt depth -2.6pt
                   \kern -6pt \intop}\nolimits_{#1}}%
        {\mathop{\vrule width 5pt height 3 pt depth -2.6pt
               \kern -6pt \intop}\nolimits_{#1}}%
         {\mathop{\vrule width 5pt height 3 pt depth -2.6pt
                  \kern -6pt \intop}\nolimits_{#1}}}

\newcommand{\Foa}{\,\,\,\text{for all }\,\,}

\begin{document}

\title{
On sphere--filling ropes 
}
\author{Henryk Gerlach, Heiko von der Mosel}
\maketitle

\begin{abstract}
What is the longest rope on the unit sphere? Intuition tells us
that the answer to this packing problem depends on the rope's
thickness. For a countably infinite number of prescribed thickness
values we construct and classify all solution curves. The simplest
ones are similar to
the seamlines of a tennis ball, others exhibit a striking resemblance
to 
 Turing patterns in chemistry, or to ordered phases
of long elastic rods stuffed into spherical shells.
\bigskip

\end{abstract}

\begin{figure}
  \includegraphics[width=4cm]{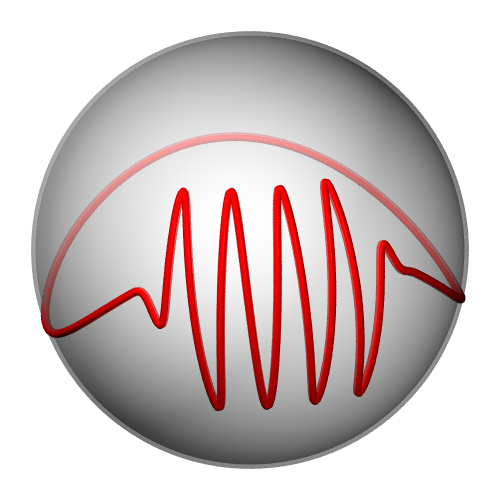}
  \hspace{1cm}\begin{minipage}[b]{7cm}
  \caption{
           Without a lower bound on the thickness 
           there is no longest curve on $\mS^2$. Inserting more and more
           oscillations into a given curve its length can be made
	   arbitrarily large.
          }\label{fig1}
  \end{minipage}

\end{figure}

\section{The problem}\label{sec:1}
What is the longest curve on the unit sphere? The most probable
answer of any mathematically inclined person to this na\"ive
question
is: There is no such
thing, since any spherical curve of finite length can be made
arbitrarily long by replacing parts of it by more and more ``wavy''
arcs; see Figure \ref{fig1}. 
Rephrasing the initial query as ``what is the longest {\it rope} on the unit
sphere?'' makes a big difference. A rope in contrast to a mathematical
curve forms a solid body with positive thickness, so that now this
question addresses a packing problem with obvious parallels in everyday
life. Is there an optimal way of winding electrical cable onto the reel?
Similarly, and economically quite relevant, can one maximize the volume of yarn
wound onto a given bobbin \cite{mayer_lenz}, or how should one store long textile fibre band most
efficiently to save storage space \cite{kyosev}?  

Common to all these packing problems, in contrast to the classic
Kepler problem of optimal sphere packing \cite{CSl88}, \cite{hales}, is
that long and slender deformable objects are to be placed into a fixed volume
or onto a given surface.  Nature displays 
fascinating packing strategies on various scales. Extremely long strands
of viral DNA are  packed very efficiently 
into the tiny phage head of bacteriophages \cite{EH77}, 
and  chromatin fibres
are folded and 
organized in various aggregates within the chromatid \cite{MuJ83}.

\begin{figure}
  \begin{center}
    \includegraphics[height=3cm]{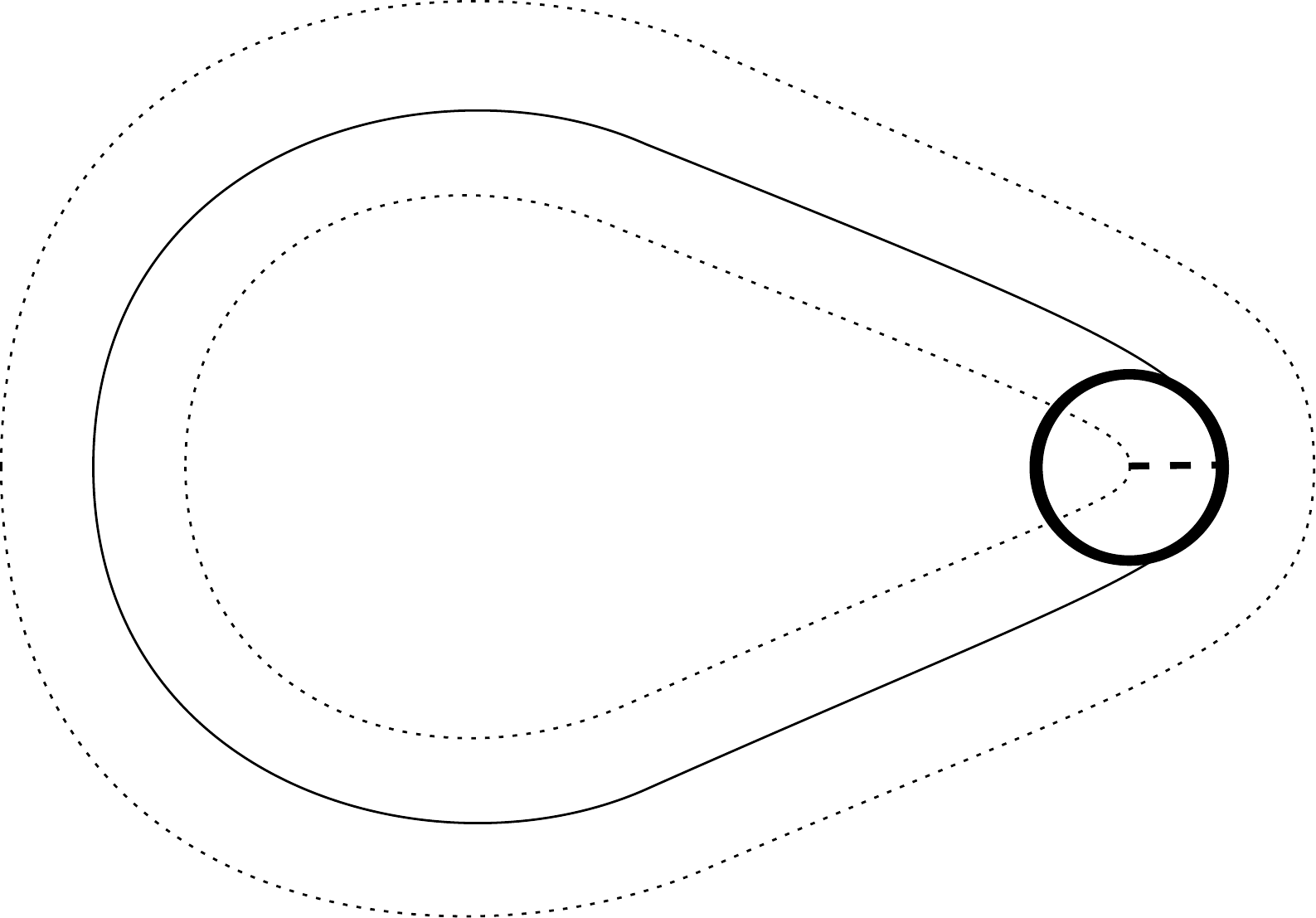} 
    \hfill
    \includegraphics[height=3cm]{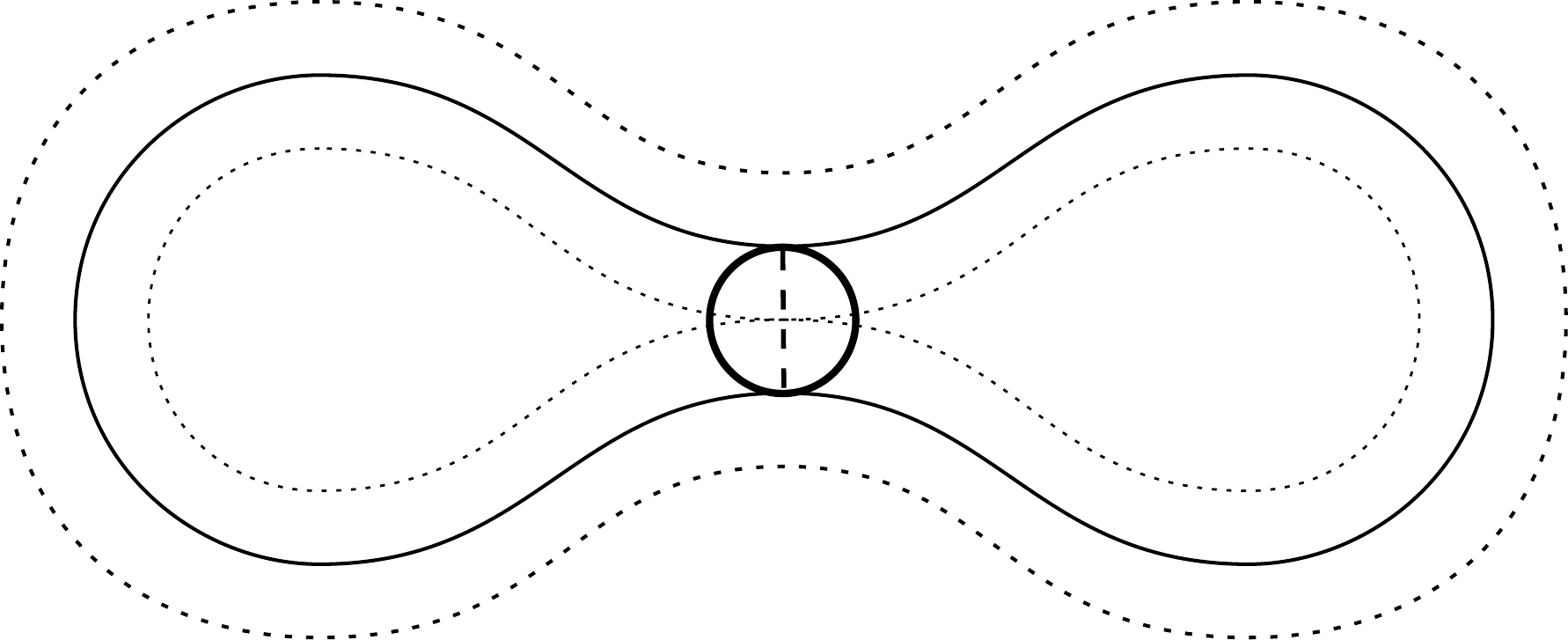}
  \end{center}
  \caption{ A positive thickness imposes a lower 
           bound both on radius of curvature (left) as well as
	   on global self-distance (right). The tubular neighbourhood 
           around the curve does not self-intersect.}
          \label{fig2}
\end{figure}

To model a rope as a mathematical curve $\gamma$ with positive 
thickness we follow the approach of Gonzalez and Maddocks \cite{GM99}
who considered all triples of distinct curve points $x,y,z\in\gamma$,
and their respective circumcircle radius $R(x,y,z)$. The smallest
of these radii determines the curve's thickness
\begin{equation}\label{1.0}
\triangle[\gamma]:=\inf_{x\not= y\not=z\not=x\atop
x,y,z\in\gamma}R(x,y,z).
\end{equation}
A positive lower bound on this quantity controls local curvature
but also prevents the curve from self-intersections; see Figure \ref{fig2}.
In fact, it equips the curve with a tubular neighbourhood of uniform
radius $\triangle[\gamma]$ without self-penetration.
It can be shown that positive thickness characterizes
the set of embedded curves with bounded curvature 
\cite[Lemmata 2 \& 3]{GMSvdM02}, \cite[Theorem 1]{SvdM03a},
and we therefore tacitly assume from now on that our curves are
simple, have positive length,  and are continuously differentiable.

With this mathematical concept of thickness at our disposal we 
can reformulate the original question of finding the longest ropes on
the unit sphere as a variational problem, where we first focus on
closed loops. 
\begin{problem*}
For a given constant $\Theta>0$ find the longest
closed curve $\gamma:\S^1\cong \R/(2\pi\Z)\to\S^2:=\{x\in\R^3:|x|=1\}$
with prescribed minimal thickness $\Theta$, i.e., with
$\triangle[\gamma]\ge \Theta.$
\end{problem*}
Before discussing the solvability of this maximization problem for
various thickness values we would like to point out that every loop
$\gamma\subset\R^3$ of positive thickness enjoys a strong 
geometric property, the presence of {\it forbidden balls}: {\it Any
open ball $B_\Theta\subset\R^3$ of radius $\Theta\le\triangle[\gamma]$
whose boundary $\partial B_\Theta$ touches the curve $\gamma$ 
tangentially in a point $p\in\gamma$, is not penetrated by the curve, that
is $B_\Theta\cap\gamma=\emptyset.$}
In fact, otherwise 
there were a point $q\in B_\Theta\cap\gamma$, and the plane spanned by the
segment $q-p$ and the tangent vector of $\gamma$ at $p$ would intersect
$B_\Theta$ in a planar disk of radius at most $\Theta.$ This disk would contain
the strictly smaller circle through $q$ and $p$ that is tangent
to the disk's
boundary in $p$. Approximating this circle by the circumcircles
of the point triples $q,p,p_i$  for some sequence $\{p_i\}$ of curve
points converging to $p$ as $i$ tends to infinity, yields a contradiction
via
$
\triangle[\gamma]\le R(p,q,p_i)<\Theta$ for sufficiently
large $i$.

A direct consequence of the presence of forbidden balls is that
Problem (P) is not solvable at all if the prescribed thickness
is strictly greater than $1$, there are simply no spherical curves
whose thickness exceeds the value $1$. Indeed,
for any point $p$ on a spherical curve $\gamma$ with thickness
$\triangle[\gamma]>1$ there exists an open ball of radius
$\triangle[\gamma]$ touching the unit sphere (and therefore $\gamma$
as well) in $p$ and containing all of the unit sphere but $p$.
However, this ball is forbidden, hence contains no curve point
so that $p$ is the only curve point on $\S^2.$ This settles Problem
(P) for $\Theta>1.$

\begin{figure}
  \begin{center}
    \includegraphics[height=4cm]{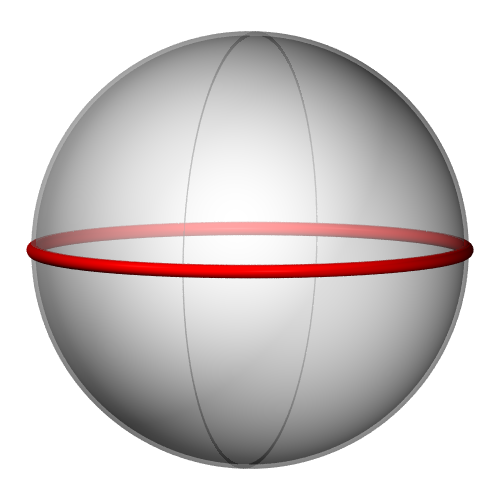}
    \hspace{1cm}
    \includegraphics[height=4cm]{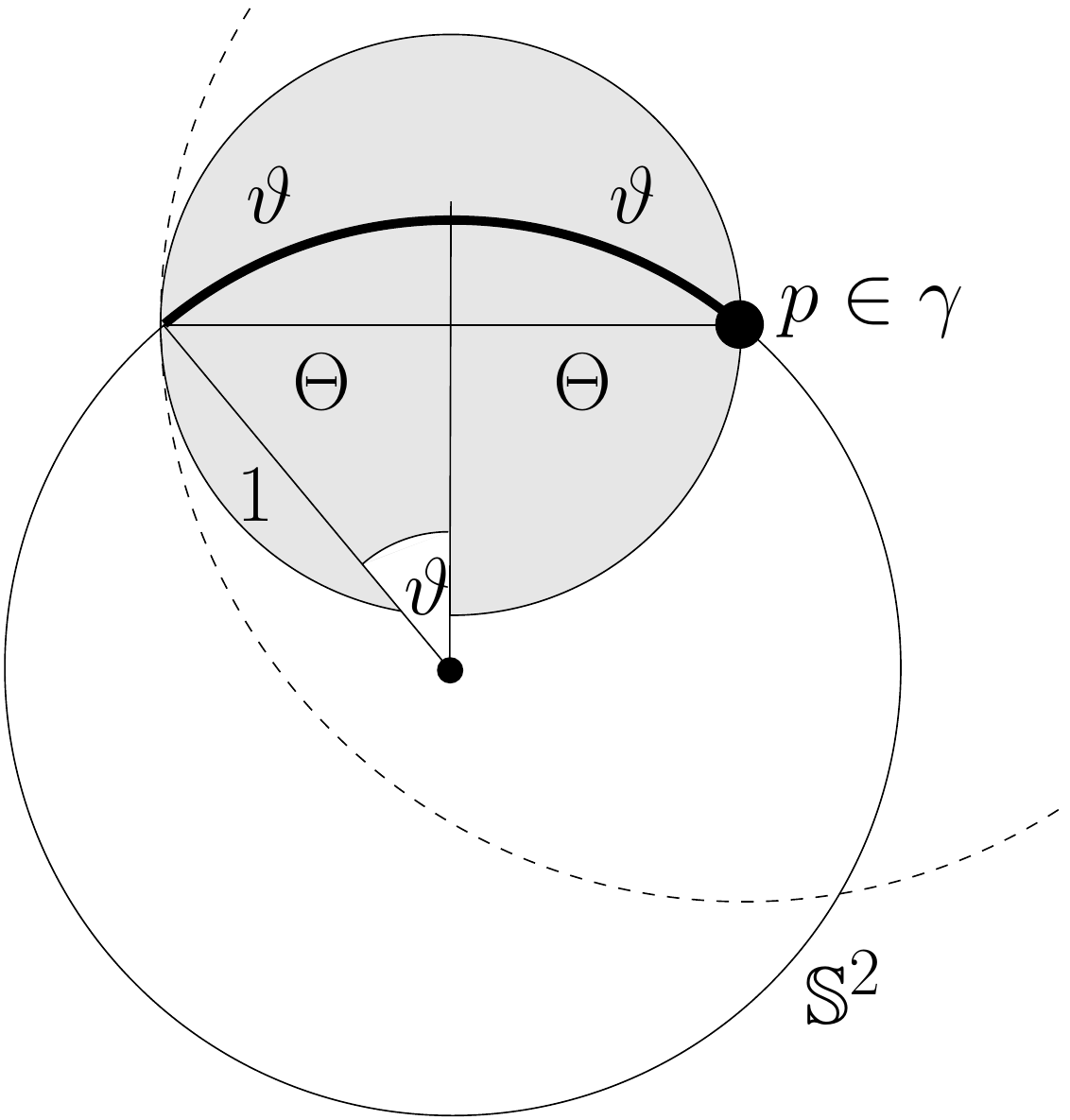}
  \end{center}
  \caption{Left: A great-circle is the thickest curve on $\mS^2$. Right:
           The unit sphere cut along a normal plane that is orthogonal to $\gamma$ at the point $p\in\gamma$. 
           The grey spatial forbidden ball 
           rotated along the dashed circle generates a forbidden 
           geodesic ball of radius $\vartheta=\arcsin\Theta$ on $\S^2$.
  }\label{fig3}
\end{figure}

If we intersect the union of all forbidden touching balls $B_\Theta$,
$\Theta\le\triangle[\gamma]$ for a loop $\gamma\subset\S^2$, with the 
unit sphere, we easily deduce (see Figure \ref{fig3}) that every curve point of a spherical
curve carries a pair of 
\begin{forbiddenballs*}
A closed spherical curve $\gamma:\S^1\to\S^2$ with (spatial)
thickness $\triangle[\gamma]\ge\Theta >0$ does not intersect
any open geodesic ball $\mathscr{B}_\vartheta(\xi):=
\{\eta\in\S^2:\distS(\eta,\xi)<\vartheta=\arcsin\Theta\}$ on $\S^2$
whose boundary $\partial\mathscr{B}_\vartheta(\xi)$
is tangent to $\gamma$ in at least one curve point.
Here $\distS(\cdot,\cdot)$ denotes the intrinsic distance on $\S^2$.
\end{forbiddenballs*}

One can imagine a bow tie consisting of two open geodesic balls
of spherical radius $\vartheta$ 
attached to the curve at their common boundary point.
This bow tie can be moved freely 
along the curve without ever  hitting any part of the curve.

The full strength of Property (FGB) is frequently used later on
to completely
classify infinitely many explicit solutions of Problem (P). For the
moment it helps us to quickly solve that problem for $\Theta=1$.
Take any point $p$ on an arbitrary spherical curve $\gamma$ with 
thickness $1$. The two forbidden open geodesic balls of spherical radius
$\vartheta=\arcsin 1=\frac{\pi}{2}$ touching $\gamma$ in $p$
are two complementary hemispheres $\S^+,\S^-$ that -- according to
(FGB) -- do not intersect $\gamma.$ Hence $\gamma$ must be the
equator as the only closed curve contained in the complement
$\S^2\setminus (\S^+\cup\S^-).$ Thus the equator is the only spherical
curve with thickness $1$ and hence -- up to congruence -- the unique
solution to Problem (P) for $\Theta=1.$

But what about other thickness values $\Theta\in
(0,1)$, is the variational problem (P) solvable at all? The 
answer is yes, and once one has analyzed the continuity properties
of the constraint $\triangle[\gamma]\ge\Theta$, this can be
proven with a direct method in the calculus of variations. The necessary
arguments for this (and for 
the constructions and classification results in Sections
\ref{sec:2} and \ref{sec:3})   are carried out in full detail in
\cite{GevdM09}. Related existence results for thick elastic rods and
ideal knots can be found in
\cite{GMSvdM02},  \cite{CKS02}, \cite{GdlL03}, \cite{GeM10}.
\begin{theorem}[Existence {\cite[Theorem 1.1]{GevdM09}}]\label{thm:1.1}
 For each prescribed minimal thickness
 $\Theta\in (0,1]$ Problem {\rm (P)} possesses (at least) one solution
 $\gamma_\Theta$. In addition, every such solution
 attains the minimal thickness, i.e., $\triangle[\gamma_\Theta]=\Theta.$
 \end{theorem}

\section{Infinitely many explicit solutions.}\label{sec:2}
Knowing that solutions exist does not necessarily mean that we know their
actual shape, unless $\Theta=1$ where we have identified  the equator as
the only solution. For general variational problems it is mostly
impossible to extract explicit information about the shape of solutions,
even uniqueness is usually a challenging issue. Here, however,
the situation is different, and this has to do with the fact that
every spherical curve $\gamma$ with positive thickness $\triangle[\gamma]
=\Theta$ carries an open tubular neighbourhood
$$
\cT_\vartheta(\gamma):=\{\xi\in\S^2:\distS(\gamma,\xi)<\vartheta=\arcsin\Theta
\},
$$ 
which equals the union of subarcs of great circles of uniform length
$2\vartheta$ on the sphere. Each of these great-arcs is centered at a 
curve point $p\in\gamma$, and is orthogonal to the
respective tangent vector of $\gamma$ at $p$. If two such great-arcs
centered at different points $p,q\in\gamma$ had a common point, then
$p$ would be contained in one of the forbidden geodesic balls touching $\gamma$
at $q$, which is excluded by Property (FGB).
Therefore this union $\cT_\vartheta(\gamma)$ of great-arcs is disjoint,
and we conclude that a curve with (spatial) thickness 
$\triangle[\gamma]=\Theta$ has the larger {\it spherical thickness}
$\vartheta=\arcsin\Theta >\Theta$.

It has been shown
more than 70 years ago by Hotelling \cite{H39}
and in more generality by  Weyl \cite{W39} that the volume of such a 
uniform tubular neighbourhood is proportional to the length $\mathscr{L}$
of its centerline. Adapted to the present situation of thick curves
on the unit sphere this classic theorem reads as
$$
\area (\mathscr{T}_\vartheta(\gamma))=2\sin\vartheta\cdot
\mathscr{L}(\gamma)
$$
for any curve $\gamma\subset\S^2$ with (spatial) thickness $\triangle[\gamma]
\ge\Theta=\sin\vartheta.$ 
Consequently, any curve $\gamma\subset\S^2$ with 
thickness $\triangle[\gamma]\ge\Theta$ whose spherical tubular neighbourhood
$\mathscr{T}_\vartheta(\gamma)$ covers {\it all} of $\S^2$, i.e., with
\begin{equation}\label{spherefillingcond}
\area(\mathscr{T}_\vartheta(\gamma))=4\pi=\area(\S^2),
\end{equation}
has {\it maximal} length among all spherical curves with 
prescribed minimal
thickness $\Theta.$ In other words, {\it sphere-filling} thick curves
provide solutions to Problem (P).

Are there any thickness values $\Theta\in (0,1)$ such that we find
sphere-filling curves of that minimal thickness, i.e., curves $\gamma
\subset\S^2$ with $\triangle[\gamma]\ge\Theta$ such that for
$\vartheta=\arcsin\Theta$ we have Relation \eqref{spherefillingcond}?

\begin{figure}
  \begin{center}
    \includegraphics[height=2.2cm]{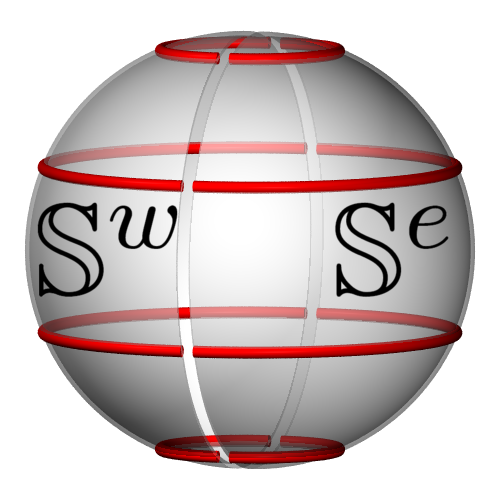}
    \includegraphics[height=2.2cm]{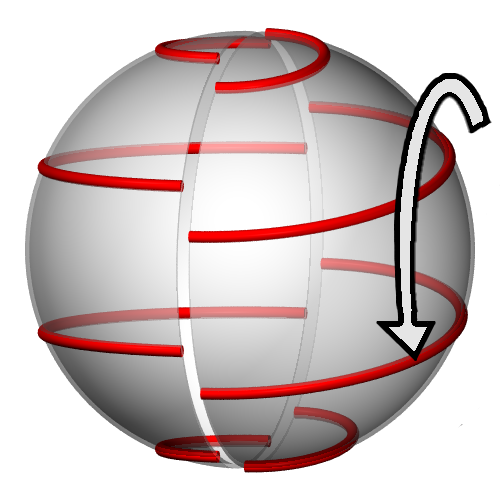}
    \includegraphics[height=2.2cm]{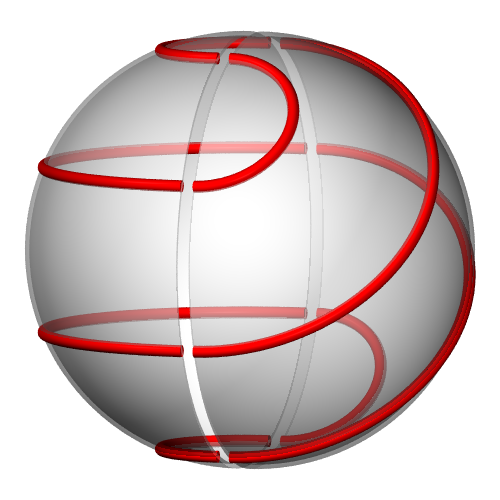}
    \includegraphics[height=2.2cm]{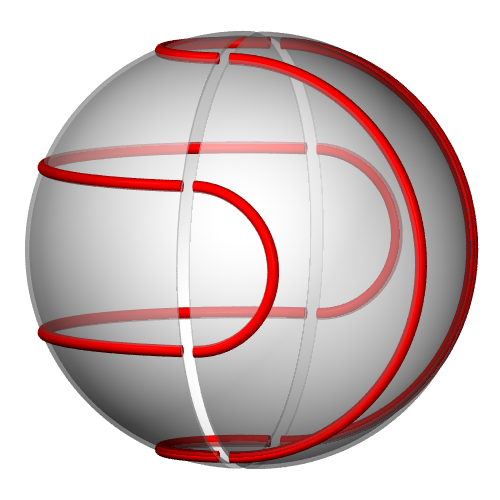}
    \includegraphics[height=2.2cm]{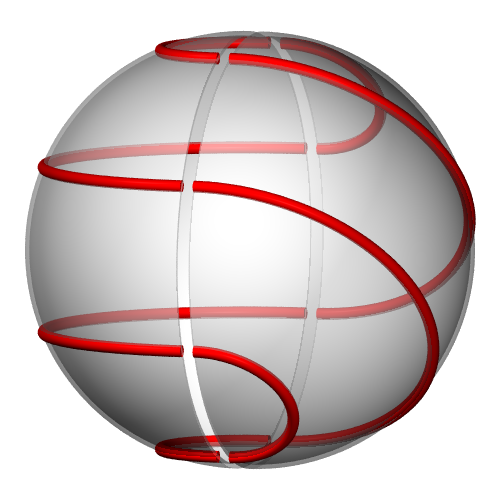}
  \end{center}
  \caption{The construction of solutions for $n=4$.
  The third and the fifth image depict the sphere-filling
  curves for turning angles $2\vartheta_4$ and $6\vartheta_4$.
  The fourth image, in contrast, contains a disconnected
  sphere-filling ensemble with two components.
  }\label{fig5}
\end{figure}

If we relax for a moment our assumption that we search for one connected
closed curve then we easily find sphere-filling ensembles of curves.
For $\Theta_n:=\sin(\pi/2n)=:\sin\vartheta_n$, $n\in\N$,
the stack of latitudinal
circles $C_i$ with $\dist_{\mS^2}(C_0,\textnormal{north pole})=\vartheta_n$
and mutual distance $\dist_{\mS^2}(C_i,C_{i-1})=2\vartheta_n$
for $i=1,\ldots,n-1,$ forms a set of $n$
spherical curves each with spherical thickness $\vartheta_n=\pi/2n$.
Their  mutually disjoint spherical tubular neighbourhoods completely
cover the sphere:
$$
\area\left[\bigcup_{i=0}^{n-1}\mathscr{T}_{\vartheta_n}(C_i)\right]=
4\pi.
$$
This collection of latitudinal circles can now be used to construct {\it
one} closed sphere-filling curve. Let us explain in detail how, for the
case $n=4$. We cut the sphere with the $4$ latitudinal circles along
a longitudinal into an eastern hemisphere $\S^e$ and a western hemisphere
$\S^w$. Each  hemisphere contains now a 
stack of $4$ latitudinal semicircles. Keeping
the western hemisphere $\S^w$ fixed we rotate the eastern hemisphere 
$\S^e$ by an angle of $2\vartheta_4=\pi/4$ such that all the
endpoints of the now turned semicircles on $\S^e$ meet endpoints of the
semicircles on $\S^w$; see Figure \ref{fig5}.

This modified collection of semicircles still has spherical thickness
$\vartheta_4=\pi/8 $ and is sphere-filling, since in the construction
the sphere-filling stack of the original latitudinal circles was only
cut orthogonally and reunited along one longitudinal, which does not change the
thickness and sphere-filling property of the ensemble. We also observe
that this new ensemble, which resembles to some extent
the seamlines on a tennis ball,
forms {\it one} closed curve, hence solves our problem -- at least for 
this particular given spatial thickness $\Theta_4=\sin\vartheta_4=
\sin(\pi/8).$  Are there other solutions for $n=4$? Why not continue rotating
the eastern hemisphere $\S^e$ against the fixed hemisphere $\S^w$
to obtain more solutions?  It turns out that a total rotation by 
$4\vartheta_4=\pi/2$ yields two connected components, which is not
what we are looking for. But turning $\S^e$ by an angle of $6\vartheta_4
=3\pi/4$ leads to another solution: a new single closed loop not
congruent to the first one; see Figure \ref{fig5}.

One can show that this procedure works well for arbitrary $n\in\N$, 
and with
a little elementary algebra\footnote{Such a construction was used
for a bead puzzle called {\it the orb} or {\it orb it} \cite{patent}
in the 1980s and the involved algebra was probably known to its
inventors.}
we can determine the exact number of
solutions:
\begin{theorem}[Explicit solutions]\label{thm:1.2}
For each $n\in\N$ and each $k\in\{0,\ldots,n-1\}$
whose greatest common divisor with $n$ equals $1$,
the construction described above starting with $n$
latitudinal 
circles $C_0,\ldots,C_{n-1}$ with
spherical distance
$$
\dist_{\mS^2}(C_0,\textnormal{north pole\,})=\vartheta_n=\frac \pi {2n},\quad
\dist_{\mS^2}(C_i,C_{i-1})=2\vartheta_n,\quad i=1,\ldots,n-1,
$$
and rotating the eastern hemisphere  against the fixed 
western hemisphere  by an angle of $2k\vartheta_n$, leads
to $\varphi(n)$ explicit piecewise circular solutions of the variational
problem {\rm (P)} for prescribed minimal thickness $\Theta_n=\sin\vartheta_n$.
\end{theorem}
Here, $\varphi$ denotes the Eulerian totient function from
number theory: $\varphi(n)$ gives the number of integers $0\le k<n$ 
so that the greatest common divisor of $k$ and $n$ equals $1$. In our example
above, $n=4$, we indeed found $\varphi(4)=2$ explicit solutions by
rotating the eastern hemisphere by the amount of $2k\vartheta_4$ for $k=1$ and for
$k=3.$

Figure \ref{fig6} depicts such sphere-filling closed curves for
various $n$, and one notices a striking resemblance with certain
so-called {\it Turing patterns} observed and analyzed in chemistry
and biology as characteristic  concentration distributions of
different substances; see, e.g., \cite{varea}. In that context,
the patterns are caused by diffusion-driven instabilities; here in
contrast, the shape of solutions is a consequence of a simple
variational principle.

\newlength\imgwidth
\setlength\imgwidth{3.4cm}
\begin{figure}
  \begin{center}
    \begin{tabular}{ccc}
      \includegraphics[height=\imgwidth]{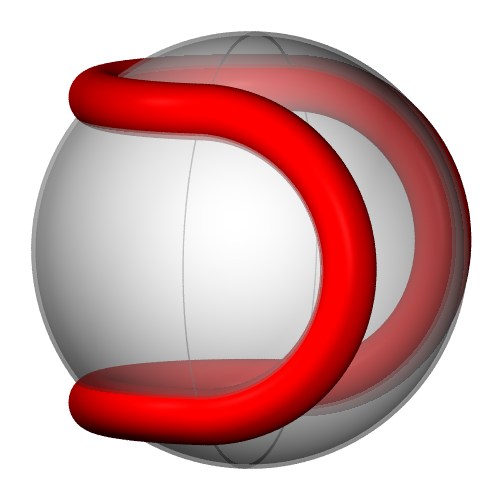} &
      \includegraphics[height=\imgwidth]{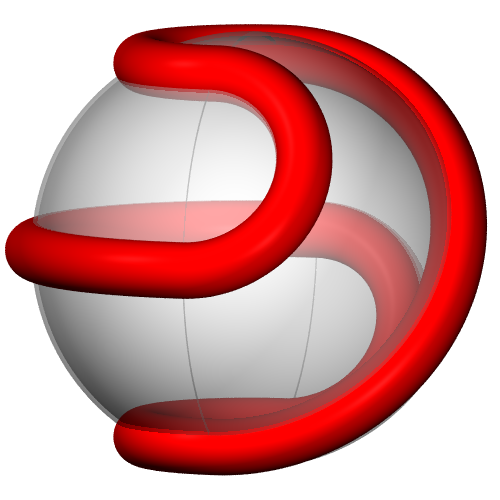} &
      \includegraphics[height=\imgwidth]{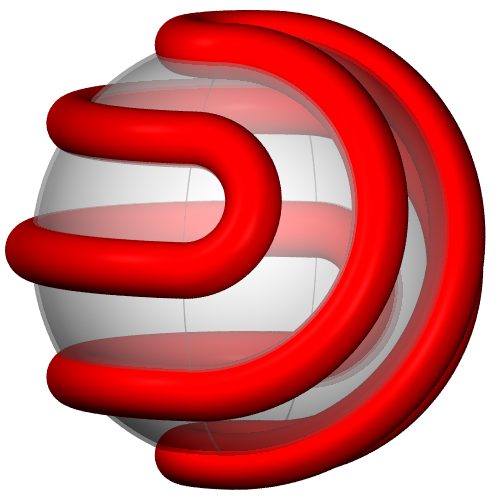} \\
      $n=2, k=1$ & 
      $n=3, k=1$ &
      $n=5, k=2$ \\
      \includegraphics[height=\imgwidth]{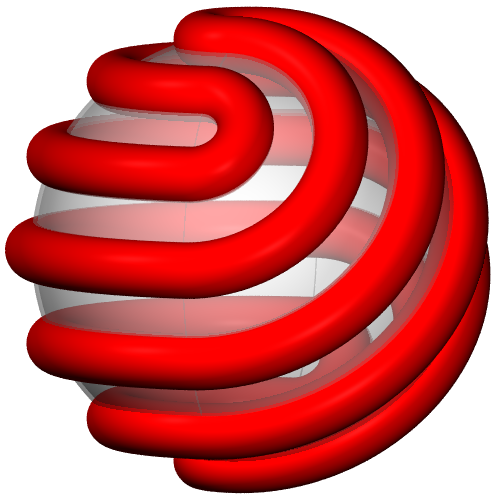} &
      \includegraphics[height=\imgwidth]{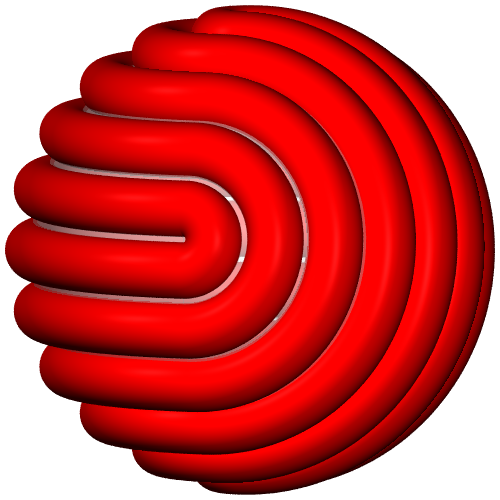} &
      \includegraphics[height=\imgwidth]{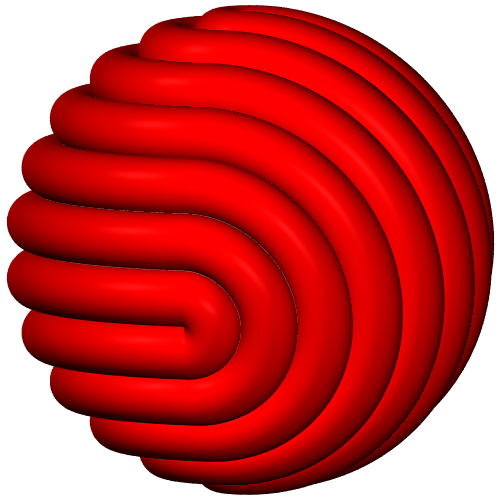} \\  
      $n=7, k=2$ &
      $n=11, k=5$ &
      $n=12, k=7$
    \end{tabular}
  \end{center}
  \caption{Various closed solutions of Problem (P). All curves are 
  visualized as tubes of a fixed radius  $\Theta=\pi/24$, which coincides
  with the actual spatial 
  thickness $\triangle[\gamma]=\Theta_n$ only 
  for the last curve with $n=12$. The remaining values of 
  spatial thickness 
  $\Theta_n$ for $n=1,\ldots,11,$ all exceed the tubes' radii
  depicted in the image. 
  }\label{fig6}
\end{figure}

\begin{figure}
  \begin{center}
     \begin{tabular}{ccc}
      \includegraphics[width=\imgwidth]{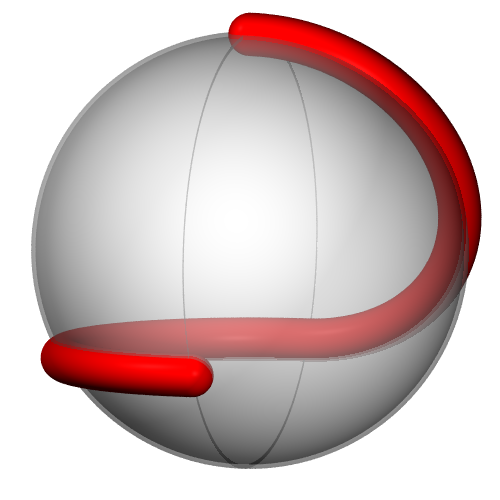} &
      \includegraphics[width=\imgwidth]{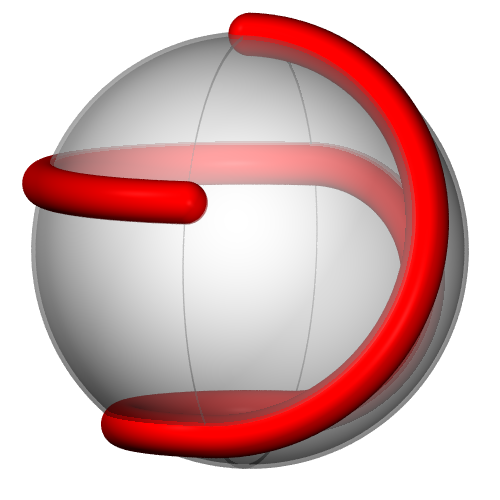} &
      \includegraphics[width=\imgwidth]{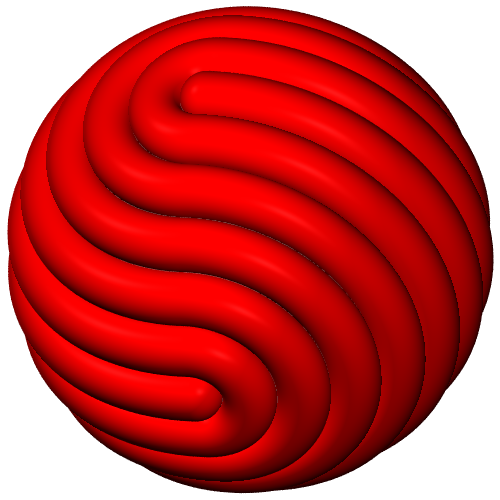} \\
       $n=3, k=1$ & $n=5, k=1$ & $n=25,k=6$  \\
      \includegraphics[width=\imgwidth]{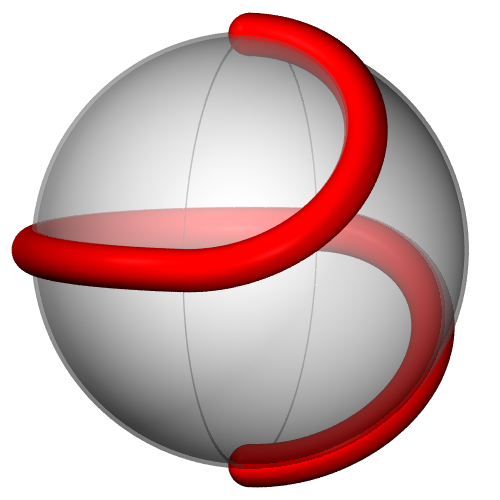} &
      \includegraphics[width=\imgwidth]{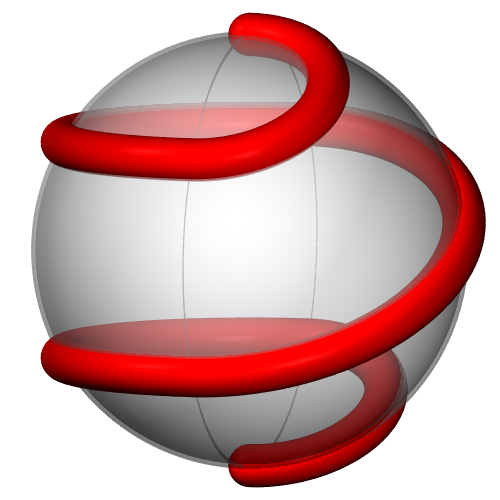} &
      \includegraphics[width=\imgwidth]{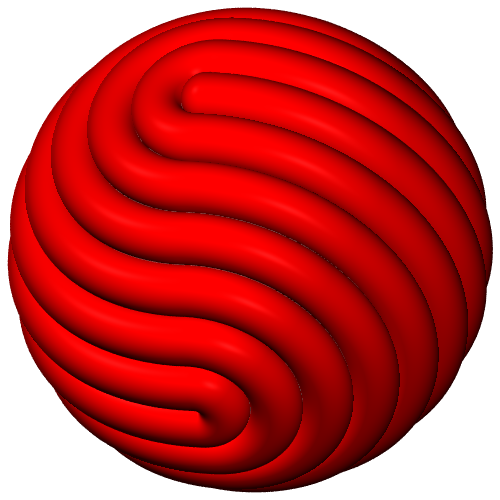} \\
       $n=4, k=0$ & $n=6, k=0$ & $n=26, k=0$  \\
    \end{tabular}
  \end{center}
  \caption{Various open solutions of Problem (P). Only the last ones 
          in each row are depicted with full spatial thickness.
          }\label{fig7}
\end{figure}

Similar constructions for thickness values $\Omega_n:=\sin(\pi/n),$
$n\in\N$, starting from an initial ensemble of semicircles together with
one or two poles on $\S^2$, lead to two disjoint families of sphere-filling
{\it open curves} distinguished by the relative position of the two endpoints
on the sphere; see Figure \ref{fig7}.
For all even $n\in\N$ the respective
open sphere-filling curves have antipodal endpoints, which is not the case
if $n$ is odd. Let us point out that these open curves occur in 
the different context of statistical physics, namely as two of three
possible configurations of ordered phases of long elastic rods densely 
stuffed into spherical shells; see  \cite{KAB06}, in particular their
figures 4a and 4c. Those studies aimed at explaining the
possible nematic order of densely packed long DNA in viral capsids.

\section{Classification of sphere-filling ropes}\label{sec:3}
For each positive integer $n$ we have constructed explicitly longest
closed ropes of thickness $\Theta_n=\sin(\pi/2n)$ 
 on the unit sphere. Are there
more? We know there are, for intermediate values $\Theta\not=\Theta_n$
by Theorem \ref{thm:1.1}, but even 
if we stick to these specific countably many values 
$\Theta_n$  of given minimal thickness we might find
more sphere-filling and thus length maximizing curves of considerably
different shapes? The answer may be surprising, but, no, up to congruence
our solutions are the only ones, and this ``uniqueness'' result is actually
a consequence of a complete classification of sphere-filling thick 
curves:
\begin{theorem}[Classification of sphere-filling loops]\label{thm:1.3}
If the spherical tubular neighbourhood $\mathscr{T}_\vartheta(\gamma),$
$\vartheta\in (0,\pi/2]$ of a closed spherical curve $\gamma\subset\S^2$
with thickness $\triangle[\gamma]\ge\Theta=\sin\vartheta$
satisfies
$$
\area(\mathscr{T}_\vartheta(\gamma))=4\pi=\area (\S^2),
$$
then there exist positive integers $n$, and $k\in\{0,\ldots,n-1\}$
with greatest common divisor equal to $1$, such that 
$\vartheta=\vartheta_n=\pi/(2n),$ $\triangle[\gamma]=\Theta_n=
\sin\vartheta_n,$ and such that $\gamma$ coincides -- up to congruence --
with one of the $\varphi(n)$ explicit solutions of Problem {\rm (P)}
exhibited in Theorem \ref{thm:1.2}.
\end{theorem}
An analogous result holds also for open curves: any sphere-filling
thick open curve must have spherical thickness $\omega_n=\arcsin\Omega_n=
\pi/n$ for some $n\in\N$, and coincides with a member of one of the
two explicitly constructed families of open spherical curves, depending
on whether $n$ is even or odd. So, if one was given the (somewhat strange)
task to produce a soccer ball of a given size
by deforming a continuous piece of thick rope 
of suitable length into an airtight 
spherical hull, then
only specific values of rope thickness are possible, 
and our theorem tells us how
one should proceed.
 There is simply no other way!

Let us explain the main ideas of the proof of this  classification result.
The presence of forbidden geodesic balls (FGB) allows us to prove
a fundamental {\it touching principle} for spherical curves $\gamma$
with positive
thickness $\triangle[\gamma]$; see Part A below.
This principle guarantees then that the number of possible local
touching situations between the curve and geodesic balls  with radius
equal to
$\triangle[\gamma]$ is very limited (Part B). The combination of these
pieces of information leads to a geometric  rigidity for
sphere-filling curves reflected in two sorts of possible 
global patterns (Part C).

\begin{figure}
  \includegraphics[width=4cm]{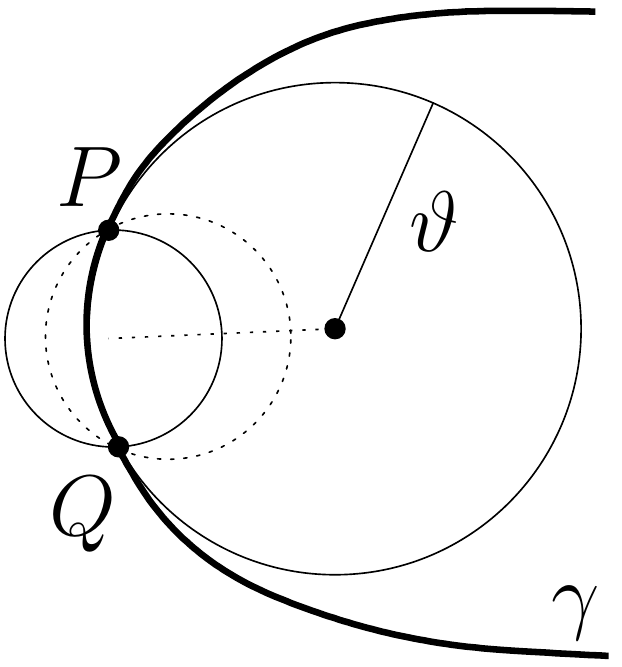}
  \hspace{1cm}\begin{minipage}[b]{7cm}
  \caption{A curve $\gamma$ with spherical thickness $\vartheta$ 
  that touches a circle of radius $\vartheta$ 
  in two non-antipodal points $P$ and $Q$, joins 
           them with a subarc of that circle. 
           Other\-wise one of the dotted circles would intersect 
	   $\gamma$ three times 
	   leading to  a lower thickness.
          }\label{fig8}
  \end{minipage}
\end{figure}

{\bf A. The touching principle}\, addresses the situation when a
spherical curve $\gamma$ with $\triangle[\gamma]\ge\Theta=
\sin\vartheta$ touches the boundary $\partial \mathscr{B}_{\vartheta}
(\xi_0)$ of a geodesic ball $\mathscr{B}_{\vartheta}(\xi_0)$ in
$\S^2$ in at least two non-antipodal points $P,Q\in\partial
\mathscr{B}_{\vartheta}
(\xi_0)$ with $\distS(P,Q)=:2\vartheta_1<2\vartheta.$ In this situation
the boundary of the strictly smaller geodesic ball $\mathscr{B}_{\vartheta_1}(\xi_1)$ for which $P$ and $Q$ are antipodal, is intersected transversally
by $\gamma$ in $P$ and $Q$, which means that the open geodesic
ball $\mathscr{B}_{\vartheta_1}(\xi_1)$ contains curve points; see Figure \ref{fig8}. 
On
the other hand, $\partial \mathscr{B}_{\vartheta_1}(\xi_1)$ contains
no further curve point $T$ different from $P$ and $Q$, since this would
imply
for the corresponding (Euclidean) circumcircle
\begin{equation}\label{star}
R(P,Q,T)\le\sin\vartheta_1 <\sin\vartheta=\Theta
\end{equation}
contradicting our assumption $\triangle[\gamma]\ge\Theta$; recall
Formula \eqref{1.0}. Consequently, there is a whole subarc of $\gamma$
connecting $P$ and $Q$ contained in $\mathscr{B}_{\vartheta_1}(\xi_1),$
but not in the original larger ball $\mathscr{B}_{\vartheta}(\xi_0)$
since this one is a forbidden ball according to (FGB).
Where can we locate this arc
within the set $\mathscr{B}_{\vartheta_1}(\xi_1)\setminus
\mathscr{B}_{\vartheta}(\xi_0)$? Sweeping out the
region $\cB_{\vartheta_1}(\xi_1)\setminus\overline{\cB_{\vartheta}(\xi_0)}
$ with intermediate geodesic circles $\partial\cB_{\vartheta_s}(\xi_s)$ with
center $\xi_s$ on the great arc connecting $\xi_0$ and $\xi_1$, and
containing $P$ and $Q$ for each $s\in [0,1]$ (so that $\vartheta_s=|P-\xi_s|$,
$\vartheta_0:=\vartheta$)
we use the same argument as the one that led to \eqref{star} to show
that there are no curve points in $\cB_{\vartheta_1}(\xi_1)\setminus
\overline{\cB_{\vartheta}(\xi_0)}.$ Thus we have proven the
\begin{touching*}
A closed spherical curve $\gamma:\S^1\to\S^2$ with (spatial)
thickness $\triangle[\gamma]\ge\Theta=\sin\vartheta$ 
that touches tangentially a geodesic
circle $\partial\cB_{\vartheta}\subset\S^2$ 
of spherical radius $\vartheta$ in two non-antipodal points $P$ and $Q$,
contains the shorter circular subarc of $\partial \cB_{\vartheta}$
connecting $P$ and $Q$.
\end{touching*}

We benefit from the touching principle since it allows us to characterize
sphere-filling curves of thickness $\Theta=\triangle[\gamma]=\sin\vartheta$
in terms of their local behaviour when touching geodesic balls of 
spherical radius $\vartheta$: For any open geodesic ball $\cB_\vartheta$
disjoint from $\gamma$ -- and there are plenty of those, e.g., all forbidden balls by (FGB) --
then one of the following three touching
situations is guaranteed for the intersection
$\mathscr{S}:=\partial \cB_\vartheta\cap
\gamma$:

\medskip

{\bf B. Possible local touching situations.}\,
\begin{enumerate}
\item[\rm (a)] $\gamma$ touches $\partial\cB_\vartheta$ in exactly
two antipodal points, i.e., $\mathscr{S}
=\{P,Q\}$ with $\distS(P,Q)=2\vartheta$, or
\item[\rm (b)] the intersection $\mathscr{S}$ 
is a relatively closed semicircle of spherical radius
$\vartheta$, or
\item[\rm (c)] this intersection $\mathscr{S}$ equals the full
geodesic circle $\partial \cB_\vartheta.$
\end{enumerate}
To see this we notice first that for a sphere-filling curve
the relatively closed intersection
$\mathscr{S}$ is nonempty, since otherwise a slightly larger
ball $\cB_{\vartheta+\epsilon}$ for some small positive 
$\epsilon$ would not contain any curve point, which leads to
$\cB_\epsilon\cap\mathscr{T}_\vartheta(\gamma)=\emptyset,$
and hence $4\pi=\area(\mathscr{T}_\vartheta(\gamma))\le\area
(\S^2\setminus\cB_\epsilon)<4\pi$  contradicting \eqref{spherefillingcond}.

Similarly, one can rule out that the set
$\mathscr{S}$ is contained in a relatively
open semicircle on $\partial\cB_\vartheta$, since then
two extremal points $\eta,\zeta\in\mathscr{S}$ realizing the diameter
of $\mathscr{S}$ would have spherical distance $\distS(\eta,\zeta)<
2\vartheta.$ This is one of the frequent occasions that the 
touching principle (TP) comes into play. It implies that the whole 
circular subarc  of $\partial\cB_\vartheta$ connecting
$\eta$ and $\zeta$, is contained in $\mathscr{S}$ and therefore
equals $\mathscr{S}$ in this situation. But the fact that $\gamma$
lifts off the geodesic circle $\partial\cB_\vartheta$ at $\eta$ and
$\zeta$ before
completing a full semicircle allows us to move the closed ball
$\overline{\cB_\vartheta}$ slightly ``away'' from $\mathscr{S}$, that is,
in the direction orthogonal to and away from the geodesic arc
connecting $\eta$ and $\zeta.$ This way we obtain a slightly shifted closed
ball of the same radius without any contact to $\gamma$, a 
situation that we have ruled out above. 
Therefore, $\mathscr{S}$ is not contained in any relatively open
semicircle on $\partial\cB_\vartheta.$

If $\mathscr{S}$ is contained in a relatively closed semicircle
we may assume that 
it contains apart from the antipodal endpoints of that semicircle
also at least one third point, otherwise we were in Situation (a)
and could stop here. Therefore, by virtue of the touching
principle (TP) $\mathscr{S}$ coincides completely with that closed
semicircle, which is Option (b). If, however, $\mathscr{S}$
is not contained in any semicircle we can simply look at one point
$q\in\mathscr{S}$ and its antipodal point $q'\in\partial\cB_\vartheta$.
If $q'$ happens to be also in $\gamma$, then both semicircles connecting
$q$ and $q'$ would contain further curve points and therefore
$\mathscr{S}=\partial\cB$ again by the touching principle,
and we end up with Option (c). 

If $q'\not\in\mathscr{S}$, then the largest
open subarc $\alpha$
of $\partial\cB_\vartheta$ containing $q'$ but no
point of $\gamma$ must be shorter than $\pi$ unless $\mathscr{S}$
is contained in a semicircle, a situation we brought to a close
before.
Applying the touching principle to the two endpoints of $\alpha$ we 
find in fact that $q'$ lies on the subarc of $\gamma$ connecting
these endpoints on $\partial\cB_\vartheta,$,  which exhausts the last
possible situation to verify that our list of situations (a)--(c)
is complete.

We are going to use the local structure established in Parts A and B to 
prove geometric rigidity of sphere-filling curves $\gamma$ with 
positive spatial thickness $\Theta=\sin\vartheta.$ 

{\bf C. Global patterns of sphere-filling curves.}\,
\begin{c1*}
If $\gamma\subset\S^2$ intersects a normal plane $E$ orthogonally in
$k$ distinct points whose mutual spherical distance equals
$2\vartheta$, then $k$ is even and $\gamma$ contains a semicircle
of spherical radius $\vartheta$ in each of the two hemispheres
bounded by $E\cap\S^2$.
\end{c1*}

\begin{c2*}
If $\gamma$ contains one latitudinal semicircle $S\subset \partial
\cB_\vartheta$, then $\vartheta=\pi/2n$ for some $n\in\N$ and
the portion of $\gamma$ in the corresponding
hemisphere consists of the whole stack of $n$ latitudinal semicircles
(including $S$) with mutual spherical distance $\vartheta.$
\end{c2*}

Before providing the proofs for these rigidity results let us
explain how we can combine these to establish the

\noindent
{\bf Proof of Theorem  \ref{thm:1.3}.}\,
The goal is to show the existence of a latitudinal semicircle
of spherical radius $\vartheta$ contained in $\gamma$ in order
to apply the global pattern (C2), which then assures that
$\gamma$ consists of a stack of latitudinal semicircles with
mutual spherical distance $\vartheta$ in one hemisphere, say in
$\S^w$. This behaviour in $\S^w$ leads to the characteristic
intersection point pattern in the longitudinal circle $\partial
\S^w$ needed in order to apply (C1), which in turn guarantees
the existence of a semicircle of radius $\vartheta$ on $\S^e$
whose endpoints meet $\partial\S^e$ orthogonally.
Again property (C2) leads to a whole stack of
$n$ equidistant semicircles now on $\S^e$.
Our construction described in Section \ref{sec:2} finally reveals
the only possible loops made of two such stacks of equidistant
latitudinal semicircles meeting $\partial\S^e=\partial\S^w$ orthogonally,
which completes the proof of the 
classification theorem. The logic of proof resembles a ride on
the merry-go-round; (C1) produces the semicircle on $\S^w$
necessary to use (C2) to obtain the stack of semicircles on $\S^w$,
which itself generates the point pattern needed to apply (C1)
on $\S^e$ and finish the task via (C2)  on $\S^e.$ The only problem is:
how do we enter the merry-go-round? We have to show 
that one portion of $\gamma$ is a semicircle
of spherical radius $\vartheta$ without assuming the intersection
point pattern needed in (C1).

Let $k$ be the integer such that
$(k-1)\vartheta <\pi\le k\vartheta.$ For a fixed point $p\in\gamma$
we walk along a unit speed geodesic ray $\eta_p$ emanating from $p$
in a direction orthogonal to $\gamma$ at $p$ in search of such 
a semicircle. The geodesic ball $\cB_\vartheta(\eta_p(\vartheta))$
is a forbidden ball by means of (FGB), i.e., 
$\gamma\cap\cB_\vartheta(\eta_p(\vartheta))=\emptyset,$
where
$\eta_p(\vartheta)$ denotes the point reached on the geodesic
ray after a spherical distance $\vartheta.$ According to
the possible local touching situations we find the desired
semicircle on $\partial\cB_{\vartheta}(\eta_p(\vartheta)\cap
\gamma$ (Option (b) or (c) in B),
unless the antipodal point $\eta_p(2\vartheta) $ is contained
in $\gamma.$ In that case we continue along the same geodesic
ray passing through $\eta_p(2\vartheta)$ orthogonally to $\gamma$,
until we either find a closed semicircle on one of the geodesic
circles $\partial\cB_\vartheta((2i-1)\vartheta))$, $i=1,\ldots,k,$
or $\vartheta=\pi/k,$ and all ``antipodal points'' $\eta_p(2i\vartheta)$
are contained in $\gamma$, so that $\eta_p(2k\vartheta)=p.$
In other words, either we have found the desired semicircle
during the walk along $\eta_p$, or we have walked once around the 
whole longitudinal circle traced out by $\eta_p$ generating $k$
equidistant points where $\gamma$ intersects $\eta_p$ orthogonally.
But this is exactly what is needed to apply (C1) to finally
establish the existence of the semicircle we are looking for.

One final comment on why this exact quantization takes place, i.e.,
why we find $\vartheta=\pi/k$ so that the walk along $\eta_p$ 
pinpointing the centers $\eta_p((2i-1)\vartheta)$ of geodesic
balls on the way, actually
leads exactly back to the starting point $p=\eta_p(2\pi).$ The 
successive localization of forbidden balls according to (FGB)
and the possible local touching situations yield the fact that
all open geodesic balls $\cB_\vartheta(\eta_p((2i-1)\vartheta)$
for $i=1,\ldots,k$, are disjoint from $\gamma.$ If, for instance,
the walk had stopped too late since the step size $\vartheta$ was
too large, $k\vartheta >\pi$,  then
$\distS(p,\eta_p(2k-1)\vartheta)=2\pi-(2k-1)\vartheta <\vartheta$,
that is, $p\in\cB_\vartheta(\eta_p(2k-1)\vartheta)\cap\gamma$, a
contradiction. A similar argument works if we had stopped our
walk too early.\qed

Let us establish the global patterns (C1) and (C2) in more
detail since they served as the core tools in the proof of our
classification theorem.

\begin{figure}
  \includegraphics[width=4cm]{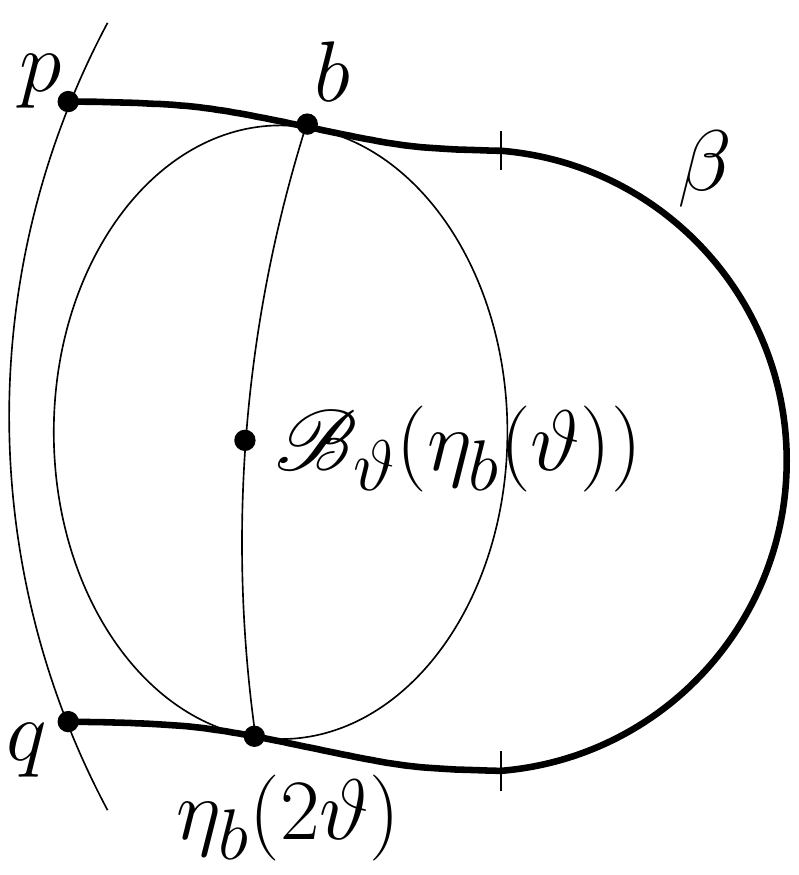}
  \hspace{1cm}\begin{minipage}[b]{7cm}
  \caption{Towards the  proof of (C1).  Geodesic rays $\eta_b$
  emanating
  from points $b$ on the subarc $\beta\subset\gamma\cap \mS^e$ to
  trace the enclosed spherical region $\mathscr{R}$. The depicted
  antipodal touching  \eqref{Anti} is excluded to hold throughout $\beta$ 
  by virtue of Brouwer's fixed point theorem.
          }\label{fig9}
  \end{minipage}
\end{figure}

We start with the 
{\bf proof of (C1)}. Here it suffices to focus
on one of the two hemispheres $\S^w$ and $\S^e$ 
determined by $E$, say on $\S^w$.
Since $\gamma$ is simple and closed the curve can leave $\S^w$ merely
as often as it enters $\S^w$, which immediately gives $k=2n$ for some
$n\in\N$. Moreover, $\S^w$ is homeomorphic to a flat disk so that
we can find nearest neighbouring exit and entrance points $p,q\in
E\cap\gamma$ with minimal spherical distance
$\distS(p,q)=2\vartheta$ such that the closed subarc $\beta\subset
\S^w\cap\gamma$ connecting $p$ and $q$ satisfies
$E\cap\beta=\{p,q\}.$
We will show that $\beta$ contains the desired semicircle
of spherical radius $\vartheta$.
Since $\gamma$ intersects $E$ orthogonally we infer from (FGB)
 that the open geodesic ball
$\cB\equiv\cB_\vartheta$ containing $p,q\in\partial\cB$ as antipodal
boundary points, is disjoint from $\gamma$.
If there were a third point $b\in\beta
\cap\partial\cB$ distinct from $p$ and $q$, then -- according
to the touching principle (TP) --
the whole semicircle on $\partial\cB$ with endpoints $p$ and $q$ would be
contained in $\gamma$, and we were done. Else we trace the
open spherical region $\mathscr{R}$ bounded by $\beta\cup
(\S^e\cap\partial\cB)$ with geodesic rays $\eta_b$  
 emanating from arbitrary points $b\in\beta$ 
orthogonally into the region $\mathscr{R}.$ Notice that $\mathscr{R}$
is disjoint from $\gamma$, and that $\eta_p(2\vartheta)
=q$ and $\eta_q(2\vartheta)=p$ where the argument of $\eta$
indicates how long one has to travel along the geodesic ray to reach
the destination point. 
In addition, the forbidden ball  property (FGB) 
implies $\gamma\cap\cB_\vartheta(
\eta_b(\vartheta))=\emptyset$ and therefore $
\cB_\vartheta(\eta_b(\vartheta))\subset\mathscr{R}$ for all points
$b\in\beta;$ see Figure \ref{fig9}.

According to Part B either
\begin{equation}\label{Anti}
\partial\cB_\vartheta(\eta_b(\vartheta))\cap\gamma =\{b,\eta_b(2\vartheta))
\quad\Foa b\in\beta,
\end{equation}
(the  antipodal situation (a)), or $\gamma$ contains a semicircle
$S_b=\partial\cB_\vartheta(\eta_b(\vartheta))\cap\gamma$ containing
itself the point $b$ for some $b\in\beta.$ This semicircle lies completely
in the western hemisphere $\S^w$,  which concludes the proof of (C1). To see
the latter assume contrariwise that there exists a point $z$ on
$\S^e\cap S_b\setminus E$. Then by connectivity also $p$ or $q$ would lie
on the semicircle $S_b$, too, which immediately implies that $S_b$ and hence
also $b\in\partial\cB\cap\beta$ lies on 
the original geodesic circle $\partial\cB$, a situation that we had excluded already.

It remains to exclude the antipodal touching \eqref{Anti}
throughout the subarc $\beta\subset\gamma\cap\S^w$. We use
Brouwer's fixed point theorem for the continuous map $f:\beta\to
\S^2$ defined by $f(b):=\eta_b(2\vartheta)$, which actually maps 
$\beta$ into $\beta$. Relation \eqref{Anti} in fact guarantees
that $f(b)\in\gamma\setminus \cB$ hence $f(b)$ is either contained
in $\S^w\cap\overline{\mathscr{R}}\cap\gamma=\beta$ in which case
we are done, or $f(b)$ lies in $\partial\cB\cap\S^e\cap\gamma.$ But 
then  the antipodal partner $b$ of $f(b)$ would also lie on
$\partial\cB$, which was excluded earlier. Consequently, Brouwer's theorem is 
applicable and leads to a fixed point $b^*=f(b^*)=\eta_{b^*}(2\vartheta)$,
which implies $2\vartheta=2\pi$ because $\eta_{b^*}$ parametrizes a unit
speed 
great circle on $\S^2$.
But this contradicts our assumption that $\vartheta\in (0,\pi/2].$
\qed

\medskip

The {\bf proof of (C2)} can be sketched as follows. Let $n$ be the
integer such that $\pi/(2n)\le\vartheta <\pi/(2n-2).$ 
According to (FGB) one has a forbidden ball:
 $\gamma\cap\cB_\vartheta=\emptyset$, and
the idea is to start from the initial semicircle $S_1:=S\subset\partial
\cB_\vartheta$ and ``scan'' the remaining part
$\S^2\setminus\overline{\cB_\vartheta}$ with unit speed 
geodesic rays $\eta_p$ 
emanating from every point $p\in S_1$ orthogonally to 
$S_1$ into the open region $\S^2\setminus\overline{\cB_\vartheta}$.
Again by (FGB) we find $\cB_\vartheta(\eta_p(\vartheta))\cap\gamma=
\emptyset$ for each such starting point $p\in S_1$. All possible 
touching situations documented in Part B guarantee the existence of at
least one second curve point $q$ on $\partial\cB_\vartheta(\eta_p(\vartheta)),$
and we claim that $q$ must be antipodal to $p$, i.e., $q=\eta_p(2\vartheta)$.
If not, then according to Options (b) or (c) in Part B, the points
$p$ and $q$ are contained in a semicircle on $\partial\cB_\vartheta(\eta_p(\vartheta))\cap\gamma.$
 But 
this semicircle is hit by neighbouring geodesic ``scanning''  rays
$\eta_r$ emanating from $r\in S_1$ for $r$ close to $p$, which would
lead to a nonempty intersection of $\gamma$ with the neighbouring
forbidden ball $\cB_\vartheta(\eta_r(\vartheta))$ contradicting (FGB).
Hence we have shown that only antipodal curve points can be generated by this
procedure: $\partial\cB_\vartheta(\eta_p(\vartheta))\cap\gamma
=\{p,\eta_p(2\vartheta))$ for all $p\in S_1$, which  produces a second
semicircle $S_2:=\bigcup_{p\in S_1}\eta_p(2\vartheta)=\partial\cB_{3\vartheta}$
contained in $\gamma$, with spherical distance $2\vartheta$ to the first
semicircle $S_1$. 

It is obvious how to continue this procedure -- now
starting the ``scanning'' rays from $S_2$ -- to obtain a whole stack
of semicircles $S_i=\partial\cB_{(2i-1)\vartheta}$ for $i=1,\ldots,n.$
If this stack is too high because the spherical thickness is too large 
with respect to $n$, i.e., if  $\pi/(2n-1)\le\vartheta <\pi/(2n-2)$, then the
stack would ``spill over'' onto the other hemisphere $\S^e$ producing
a final semicircle $S_n$ contained in $\S^e\cap\gamma$ 
with spherical radius
$(2n-1)\vartheta-\pi\in [0,\vartheta)$ that is too small: it
contradicts the spherical thickness $\vartheta=\arcsin\triangle[\gamma]$
of $\gamma$ since its curvature is to large. 
If the stack is not high enough ($\pi/2n <\vartheta <
\pi/(2n-1)$) then the last semicircle $S_n$ is still on the correct
hemisphere $\S^w$ but has spherical radius $\pi-(2n-1)\vartheta\in
(0,\vartheta)$, which is again to small for the thick curve $\gamma.$\qed

\section{Final Remark}\label{sec:4}
For an infinite countable number of thickness values we have established
a complete picture of the solution set for Problem (P) using
the sphere-filling property to a large extent. The general existence
theorem, Theorem \ref{thm:1.1}, however, guarantees the 
existence of longest ropes on the unit sphere also for all
intermediate thickness values $\Theta\not=\Theta_n$. What are 
their actual shapes? Theorem \ref{thm:1.3} ascertains that those
solutions cannot be sphere-filling. In \cite{GevdM09} we constructed
a comparison curve that could serve as a promising candidate
for prescribed minimal thickness $\Theta\in (\Theta_2,\Theta_1)$,
but this question remains to be investigated,
as well as the interesting connections to Turing patterns and the
statistical behaviour of long elastic rods under spherical confinement
mentioned in Sections 1 and 2. In addition, if one substitutes 
the unit sphere by other supporting manifolds such as the standard
torus, then the issue of analyzing the shapes of optimally packed
ropes is wide open.

{\scriptsize

\large
\vspace{1cm}
 \begin{minipage}{55mm}
    {\sc Henryk Gerlach}\\
     IMB / LCVMM \\
Station 8 \\
Facult\'e des Sciences de Base \\
EPFL Lausanne\\
CH-1015 Lausanne\\
  SWITZERLAND\\
    E-mail: {\tt henryk.gerlach@}\\{\tt epfl.ch}\\
 \end{minipage}\qquad
   \begin{minipage}{50mm}
     {\sc Heiko von der Mosel}\\
             Institut f\"ur Mathematik\\
           RWTH Aachen University\\ 
     	   Templergraben 55\\ 
   	   D-52062 Aachen\\ 
   	   GERMANY\\
     E-mail: {\tt heiko@}\\{\tt instmath.rwth-aachen.de}\\
   				       \end{minipage}

}

\begin{thebibliography}{99}

\bibitem{CKS02}
Cantarella, J.; Kusner, R.B.; Sullivan, J.M.
On the minimum ropelength of knots and links.
Inv. math. {\bf 150} (2002), 257--286.

\bibitem{CSl88}
Conway, J.H.; Sloane, N.J.A. {\it Sphere Packings, Lattices and Groups.}
Springer, New York 1988.

\bibitem{EH77} Earnshaw, W.C.; Harrison, S.C.
DNA arrangement in isometric phage heads.
Nature {\bf 268} (1977), 598--602.

\bibitem{GeM10} Gerlach, H.;   Maddocks, J.H. 
{Existence of ideal knots in $\mS^3$}.
	in preparation

\bibitem{GevdM09} Gerlach, H.; von der Mosel,
H. What are the longest ropes on the unit sphere?
Preprint Nr. 32 Inst. f. Mathematik, RWTH Aachen University 2009.

\bibitem{GdlL03} Gonzalez, O.;   de la Llave, R.
       {Existence of Ideal Knots}.
        J. Knot Theory  Ramifications {\bf 12} (2003), 123--133.

\bibitem{GM99} Gonzalez, O.;   Maddocks, J.H.
{Global Curvature, Thickness and the Ideal Shapes of Knots}. 
	Proc. Natl. Acad.  Sci. USA {\bf 96} (1999), 4769--4773.

\bibitem{GMSvdM02} Gonzalez, O.;  Maddocks, J.H.;  Schuricht,
F.; von der Mosel, H.
       {Global curvature and self-contact of nonlinearly elastic curves and rods}.
         Calc. Var. {\bf 14} (2002), 29--68.

\bibitem{hales}
Hales, T.C. Cannonballs and honeycombs. Notices AMS {\bf 47}
(2000), 440--449.

\bibitem{H39} Hotelling, H. {Tubes and Spheres in $n$-Spaces}. Amer. J. 
Math. {\bf  61} (1939), 440--460.

\bibitem{KAB06}
Katzav, E.; Adda-Bedia, M.; Boudaoud, A. 
{A statistical approach to close packing of elastic rods and to DNA packaging in viral capsids}.
Proceedings of the National Academy of Sciences, USA,
  {\bf 103} (2006), 18900--18904.

\bibitem{kyosev}
Kyosev, Y. Numerical analysis for sliver winding process with
additional can motion. In: {\it 5th International Conference
Textile Science 2003 TEXSCI 2003}, ISBN 80-7083-711-X, TU-Liberec, 
Czech Republic (2003), pp. 330-334.

\bibitem{mayer_lenz}
Mayer, M.; Lenz, F. Method and apparatus for winding a yarn into a
package. US Patent 6186435 (issued 2001).

\bibitem{MuJ83}
Mullinger, A.M.; Johnson, R.T.
Units of chromosome replication and packing.
J. Cell Sci. {\bf 64} (1983), 179--193.

\bibitem{SvdM03a} Schuricht, F.; von der Mosel, H.
{Global curvature for rectifiable loops}.
        Math. Z. {\bf 243} (2003), 37--77.

\bibitem{varea} Varea, C.; Aragon, J.L.; Barrio, R.A.
Turing patterns on a sphere.
Phys. Rev. E {\bf 60} (1999), 4588--4592.

\bibitem{W39} Weyl, H.  {On the volume of tubes}. Amer.J. Math.
{\bf 61} (1939), 461--472.

\bibitem{patent} Wiggs, C.C.; Taylor, C.J.C. {Bead puzzle}. US Patent D269629 (issued 1983).
	
\end{thebibliography}
\end{document}